\documentclass[english]{article}
\usepackage[T1]{fontenc}
\usepackage[latin9]{inputenc}
\usepackage{amssymb}

\makeatletter
\newcommand{\lyxaddress}[1]{
\par {\raggedright #1
\vspace{1.4em}
\noindent\par}
}

\makeatother

\usepackage{babel}
\begin{document}

\title{Additional congruence conditions on the number of terms in sums of
consecutive squared integers equal to squared integers }

\author{Vladimir Pletser}

\maketitle

\lyxaddress{European Space Research and Technology Centre, ESA-ESTEC P.O. Box
299, NL-2200 AG Noordwijk, The Netherlands; E-mail: Vladimir.Pletser@esa.int}
\begin{abstract}
\noindent The problem of finding all the integer solutions in $a,M$
and $s$ of sums of $M$ consecutive integer squares starting at $a^{2}\geq1$
equal to squared integers $s^{2}$, has no solutions if $M\equiv3,5,6,7,8$
or $10\left(mod\,12\right)$ and has integer solutions if $M\equiv0,9,24$
or $33\left(mod\,72\right)$; or $M\equiv1,2$ or $16\left(mod\,24\right)$;
or $M\equiv11\left(mod\,12\right)$. In this paper, additional congruence
conditions are demonstrated on the allowed values of $M$ that yield
solutions to the problem by using Beeckmans' eight necessary conditions,
refining further the possible values of $M$ for which the sums of
$M$ consecutive integer squares equal integer squares.
\end{abstract}


\section{Introduction}

The general problem of sums of $M$ consecutive integer squares starting
from $a^{2}\geq1$ being equal to integer squares $s^{2}$ involves
solving a single Diophantine quadratic equation in three variables
$M,a$ and $s$ that reads

\begin{equation}
\sum_{i=0}^{M-1}\left(a+i\right)^{2}=M\left[\left(a+\frac{M-1}{2}\right)^{2}+\frac{M^{2}-1}{12}\right]=s^{2}\label{eq:53-3}
\end{equation}

where $M>1,a,s\in\mathbb{Z}^{+}$, $i\in\mathbb{Z}^{*}$. 

With the notations of this paper, Alfred investigated \cite{key-1-22}
several necessary conditions using basic congruence equations of $M$.
Philipp \cite{key-1-23} extended Alfred's work to confirm some of
his findings. Beeckmans showed \cite{key-1-4} that with eight necessary
conditions on $M$, all values of $M$ could be found, given in a
Table with values of $M<1000$, and noted two cases for $M=25$ and
$842$ that fulfilled the eight necessary conditions but did not solve
the problem. \medskip{}

The eight necessary conditions given by Beeckmans \cite{key-1-4}
on the value of $M$ for (\ref{eq:53-3}) to hold can be summarized
as follows, with $e,\alpha\in\mathbb{Z}^{+}$:

(C1.1) If $M\equiv0\left(mod\,2^{e}\right)$, then $e\equiv1\left(mod\,2\right)$.

(C1.2) If $M\equiv0\left(mod\,3^{e}\right)$, then $e\equiv1\left(mod\,2\right)$.

(C1.3) If $M\equiv-1\left(mod\,3^{e}\right)$, then $e\equiv1\left(mod\,2\right)$.

(C2) If $p>3$ is prime, $M\equiv0\left(mod\, p^{e}\right)$, $e\equiv1\left(mod\,2\right)$,
then $p\equiv\pm1\left(mod\,12\right)$.

(C3) If $p\equiv3\left(mod\,4\right),p>3$ is prime, $M\equiv-1\left(mod\, p^{e}\right)$,
then $e\equiv0\left(mod\,2\right)$.

(C4.1) $M\neq3\left(mod\,9\right)$.

(C4.1) $\forall\alpha\geq2$, $M\neq\left(2^{\alpha}-1\right)\left(mod\,2^{\alpha+2}\right)$.

(C4.1) $\forall\alpha\geq2$, $M\neq2^{\alpha}\left(mod\,2^{\alpha+2}\right)$.

\medskip{}

Extending Beeckmans' work, it was demonstrated \cite{key-1} using
these conditions that, for (\ref{eq:53-3}) to hold, $M$ cannot be
congruent to $3,5,6,7,8$ or $10\left(mod\,12\right)$. 

It was further demonstrated \cite{key-1}, independently from Beeckmans'
conditions, that for (\ref{eq:53-3}) to hold, $M\equiv0,1,2,4,9$
or $11\left(mod\,12\right)$. Furthermore,

- if $M\equiv0\left(mod\,12\right)$, then $M\equiv0$ or $24\left(mod\,72\right)$; 

- if $M\equiv1\left(mod\,12\right)$, then $M\equiv1\left(mod\,24\right)$; 

- if $M\equiv2\left(mod\,12\right)$, then $M\equiv2\left(mod\,24\right)$; 

- if $M\equiv4\left(mod\,12\right)$, then $M\equiv16\left(mod\,24\right)$;

- if $M\equiv9\left(mod\,12\right)$, then $M\equiv9$ or $33\left(mod\,72\right)$.

These are called allowed values of $M$. The values of $M$ yielding
solutions to (\ref{eq:53-3}) are given in \cite{key-1-14}.

In this paper, additional congruence conditions on the allowed values
of $M$ such that (\ref{eq:53-3}) holds are demonstrated by using
Beeckmans' necessary conditions, refining further the possible values
of $M$ for which the sums of $M$ consecutive integer squares equal
integer squares.

In this paper, the notation $A\left(mod\, B\right)\equiv C$ is equivalent
to $A\equiv C\left(mod\, B\right)$ and $A\equiv C\left(mod\, B\right)\Rightarrow A=Bk+C$
means that, if $A\equiv C\left(mod\, B\right)$, then $\exists k\in\mathbb{Z}^{+}$
such that $A=Bk+C$. By convention, $\sum_{j=inf}^{sup}f\left(j\right)=0$
if $sup<inf$ .

\section{Additional congruence conditions on $M$}

Additional congruence necessary conditions can be found using Beeckmans'
eight necessary conditions written in all generality as follows. 

For $M>1,m,m_{1},A,B,p,e,q,\alpha,j\in\mathbb{Z}^{+}$, $i,\mu\in\mathbb{Z}^{*}$,
$0\leq\mu\leq11$ and noting $M=12m+\mu=12Am_{1}+B\mu$, one has $\forall i\geq0$
(unless indicated otherwise):

(C1.1) $M\left(\right)\neq0mod\,2^{e}$ $\Rightarrow$ $M\neq2^{e}q$
with $e\equiv0\left(mod\,2\right)$ $\Rightarrow$ $e=2i$, $q\equiv1\left(mod\,2\right)$ 

\quad{}\qquad{}$\Rightarrow$ $m_{1}\neq\left(2^{2i}q-B\mu\right)/12A$,
$\forall i>1$.

(C1.2) $M\neq0\left(mod\,3^{e}\right)$ $\Rightarrow$ $M\neq3^{e}q$
with $e\equiv0\left(mod\,2\right)$ $\Rightarrow$ $e=2i$, $q\neq0\left(mod\,3\right)$ 

\quad{}\qquad{}$\Rightarrow$ $m_{1}\neq\left(3^{2i}q-B\mu\right)/12A$,
$\forall i>0$.

(C1.3) $M\neq-1\left(mod\,3^{e}\right)$ $\Rightarrow$ $M\neq3^{e}q-1$
with $e\equiv0\left(mod\,2\right)$ $\Rightarrow$ $e=2i$, 

\quad{}\qquad{}$q\neq0\left(mod\,3\right)$ $\Rightarrow$ $m_{1}\neq\left(3^{2i}q-B\mu-1\right)/12A$,
$\forall i>0$.

(C2) If $p>3$ is prime, then $M\neq p^{e}q$ $\Rightarrow$ $m_{1}\neq\left(p^{e}q-B\mu\right)/12A$
with either

\quad{}\qquad{}if $p\equiv\pm1\left(mod\,12\right)$, $e\equiv0\left(mod\,2\right)$
$\Rightarrow$ $e=2i$, $\forall i>0$, or

\quad{}\qquad{}if $p\neq\pm1\left(mod\,12\right)$, $e\equiv1\left(mod\,2\right)$
$\Rightarrow$ $e=2i+1$.

(C3) If $p\equiv3\left(mod\,4\right)$, $p>3$ is prime and $e\equiv1\left(mod\,2\right)$
$\Rightarrow$ $e=2i+1$

\enskip{}\qquad{}$\Rightarrow$ $M+1\neq p^{2i+1}q$ $\Rightarrow$
$m_{1}\neq\left(p^{2i+1}q-B\mu-1\right)/12A$.

(C4.1) $M\neq3\left(mod\,9\right)$ $\Rightarrow$ $M\neq9q+3$ $\Rightarrow$
$m_{1}\neq\left(9q-B\mu+3\right)/12A$.

(C4.2) $M\neq\left(2^{\alpha}-1\right)\left(mod\,2^{\alpha+2}\right)$
$\Rightarrow$ $M\neq2^{\alpha}\left(4q+1\right)-1$

\quad{}\qquad{}$\Rightarrow$ $m_{1}\neq\left(2^{\alpha}\left(4q+1\right)-B\mu-1\right)/12A$,
$\forall\alpha\geq2$. 

(C4.3) $M\neq2^{\alpha}\left(mod\,2^{\alpha+2}\right)$ $\Rightarrow$
$M\neq2^{\alpha}\left(4q+1\right)$ $\Rightarrow$ $m_{1}\neq\left(2^{\alpha}\left(4q+1\right)-B\mu\right)/12A$, 

\quad{}\qquad{}$\forall\alpha\geq2$.

\medskip{}

These necessary conditions are applied to $M=12m+\mu$ for each case
of $\mu=0,1,2,4,9,11$. For brevity, conditions that are not applicable,
i.e. with which $M$ is always compliant, are not indicated. 

For $\mu=0$, $A=2$ and taking in all generality $M\equiv0\left(mod\,24\right)$
$\Rightarrow$ $M=24m_{1}$:

(C1.1) $q\equiv3\left(mod\,6\right)$, $m_{1}\neq2^{2i-3}\left(mod\,2^{2i-2}\right)$,
$\forall i>1$.

(C1.2) if $q\equiv1\left(mod\,3\right)$ $\Rightarrow$ $q\equiv16\left(mod\,24\right)$,
$m_{1}\neq\left(2\times3^{2i-1}\right)\left(mod\,3^{2i}\right)$, 

\quad{}\qquad{}if $q\equiv2\left(mod\,3\right)$ $\Rightarrow$
$q\equiv8\left(mod\,24\right)$, $m_{1}\neq3^{2i-1}\left(mod\,3^{2i}\right)$.

(C2) $q\equiv0\left(mod\,24\right)$, $m_{1}\neq0\left(mod\, p^{e}\right)$
with either

\quad{}\qquad{}if $p\equiv\pm1\left(mod\,12\right)$, $e\equiv0\left(mod\,2\right)$
$\Rightarrow$ $e=2i>0$, $\forall i>0$, or

\quad{}\qquad{}if $p\neq\pm1\left(mod\,12\right)$, $e\equiv1\left(mod\,2\right)$
$\Rightarrow$ $e=2i+1$.

(C3) $q\equiv p\left(mod\,24\right)$, $m_{1}\neq\left[\left(p^{2i+2}-1\right)/24\right]\left(mod\, p^{2i+1}\right)$.

(C4.1) $q\equiv5\left(mod\,8\right)$, $m_{1}\neq2\left(mod\,3\right)$.

(C4.3) $q\equiv2\left(mod\,3\right)$, $m_{1}\neq\left(3\times2^{\alpha-3}\right)\left(mod\,2^{\alpha-1}\right)$,
$\forall\alpha\geq3$.

For $\mu=1$, $A=2$, $B=1$, yielding $M\equiv1\left(mod\,24\right)$
$\Rightarrow$ $M=24m_{1}+1$:

(C2) if $p\equiv\pm1\left(mod\,12\right)$ $\Rightarrow$ $q\equiv1\left(mod\,24\right)$, 

\quad{}\qquad{}$m_{1}\neq\left[\left(p^{2i}-1\right)/24\right]\left(mod\, p^{2i}\right)$, 

\quad{}\quad{}if $p\neq\pm1\left(mod\,12\right)$ $\Rightarrow$
$q\equiv p\left(mod\,24\right)$, 

\quad{}\qquad{}$m_{1}\neq\left[\left(p^{2i+2}-1\right)/24\right]\left(mod\, p^{2i+1}\right)$
(%
\footnote{\noindent Beeckmans notes \cite{key-1-4} that $M=25$ cannot be rejected
using his necessary conditions. It is found here that for $p=5$ and
$i=0$, $m\neq1$ and $M=25$ can be rejected.%
}).

(C3) $q\equiv2p\left(mod\,24\right)$, $m_{1}\neq\left[\left(p^{2i+2}-1\right)/12\right]\left(mod\, p^{2i+1}\right)$.

\medskip{}

For $\mu=2$, $A=2$, $B=1$, yielding $M\equiv2\left(mod\,24\right)$
$\Rightarrow$ $M=24m_{1}+2$:

(C1.3) if $q\equiv1\left(mod\,3\right)$ $\Rightarrow$ $q\equiv19\left(mod\,24\right)$,
$m_{1}\neq\left[\left(19\times3^{2i-1}-1\right)/8\right]\left(mod\,3^{2i}\right)$, 

\quad{}\qquad{}if $q\equiv2\left(mod\,3\right)$ $\Rightarrow$
$q\equiv11\left(mod\,24\right)$, $m_{1}\neq\left[\left(11\times3^{2i-1}-1\right)/8\right]\left(mod\,3^{2i}\right)$.

(C2) if $p\equiv\pm1\left(mod\,12\right)$ $\Rightarrow$ $q\equiv2\left(mod\,24\right)$,
$m_{1}\neq\left[\left(p^{2i}-1\right)/12\right]\left(mod\, p^{2i}\right)$,

\quad{}\quad{}$\forall i>0$, 

\quad{}\quad{}if $p\neq\pm1\left(mod\,12\right)$ $\Rightarrow$
$q\equiv2p\left(mod\,24\right)$, $m_{1}\neq\left[\left(p^{2i+2}-1\right)/12\right]\left(mod\, p^{2i+1}\right)$.

(C3) $q\equiv3p\left(mod\,24\right)$, $m_{1}\neq\left[\left(p^{2i+2}-1\right)/8\right]\left(mod\, p^{2i+1}\right)$.

\medskip{}

For $\mu=4$, $A=2$, $B=4$, yielding $M\equiv16\left(mod\,24\right)$
$\Rightarrow$ $M=24m_{1}+16$:

(C1.1) $q\equiv1\left(mod\,6\right)$, $m_{1}\neq\left[2\left(2^{2i-4}-1\right)/3\right]\left(mod\,2^{2i-2}\right)$
(%
\footnote{The sequence $\left(2^{2n}-1\right)/3=1,5,21,85,341,1365,...$ is
given in \cite{key-2}%
}) or

\quad{}\qquad{}$m_{1}\neq\left(2\sum_{j=3}^{i}2^{2\left(j-3\right)}\right)\left(mod\,2^{2i-2}\right)$
, $\forall i>1$. 

(C2) if $p\equiv\pm1\left(mod\,12\right)$ $\Rightarrow$ $q\equiv32\left(mod\,48\right)$,
$m_{1}\neq\left[2\left(p^{2i}-1\right)/3\right]\left(mod\, p^{2i}\right)$,

\hspace{0.75cm}if $p\neq\pm1\left(mod\,12\right)$ $\Rightarrow$
$q\equiv16\left(mod\,48\right)$, $m_{1}\neq\left[2\left(p^{2i+1}-1\right)/3\right]\left(mod\, p^{2i+1}\right)$.

(C3) $q\equiv17p\left(mod\,24\right)$, $m_{1}\neq\left[17\left(p^{2i+2}-1\right)/24\right]\left(mod\, p^{2i+1}\right)$.

(C4.3) if $\alpha\equiv0\left(mod\,2\right)$ $\Rightarrow$ $q\equiv0\left(mod\,3\right)$,

\quad{}\qquad{}$m_{1}\neq\left[2\left(2^{\alpha-4}-1\right)/3\right]\left(mod\,2^{\alpha-1}\right)$
or $m_{1}\neq\left(2\sum_{j=3}^{\alpha/2}2^{2\left(j-3\right)}\right)\left(mod\,2^{\alpha-1}\right)$,

\quad{}\qquad{}if $\alpha\equiv1\left(mod\,2\right)$ $\Rightarrow$
$q\equiv1\left(mod\,3\right)$, 

\quad{}\qquad{}$m_{1}\neq\left[2\left(5\times2^{\alpha-4}-1\right)/3\right]\left(mod\,2^{\alpha-1}\right)$
(%
\footnote{The sequence $\left(5\times2^{2n}-2\right)/3=1,6,26,106,426,1706,...$
is given in \cite{key-3}.%
}) or

\quad{}\qquad{}$m_{1}\neq\left(2^{\alpha-3}+\sum_{j=2}^{\left(\alpha-1\right)/2}2^{2j-3}\right)$,
$\forall\alpha\geq3$. 

\medskip{}

For $\mu=9$, $A=2$, $B=1$ and taking in all generality $M\equiv9\left(mod\,24\right)$
$\Rightarrow$ $M=24m_{1}+9$:

(C1.2) if $q\equiv1\left(mod\,3\right)$ $\Rightarrow$ $q\equiv1\left(mod\,24\right)$,
$m_{1}\neq\left[3\left(3^{2i-2}-1\right)/8\right]\left(mod\,3^{2i}\right)$
(%
\footnote{The sequence $\left(3^{2n}-1\right)/8=1,10,91,820,7381,...$ is given
in \cite{key-2-5}.%
}), 

\qquad{}\qquad{}or $m_{1}\neq\left(\sum_{j=0}^{i-2}3^{2j}\right)\left(mod\,3^{2i}\right)$,

\quad{}\qquad{}if $q\equiv2\left(mod\,3\right)$ $\Rightarrow$
$q\equiv17\left(mod\,24\right)$, 

\qquad{}\qquad{}$m_{1}\neq\left(2\times3^{2i-1}+\left[3\left(3^{2i-2}-1\right)/8\right]\right)\left(mod\,3^{2i}\right)$.

(C2) if $p\equiv\pm1\left(mod\,12\right)$ $\Rightarrow$ $q\equiv9\left(mod\,24\right)$,
$m_{1}\neq\left[3\left(p^{2i}-1\right)/8\right]\left(mod\, p^{2i}\right)$;

\hspace{0.75cm}if $p\neq\pm1\left(mod\,12\right)$, 

\quad{}\qquad{}if $p\equiv5\left(mod\,24\right)$ $\Rightarrow$
$q\equiv21\left(mod\,24\right)$, $m_{1}\neq\left[\left(7p^{2i+1}-3\right)/8\right]\left(mod\, p^{2i+1}\right)$,

\quad{}\qquad{}if $p\equiv7\left(mod\,24\right)$ $\Rightarrow$
$q\equiv15\left(mod\,24\right)$, $m_{1}\neq\left[\left(5p^{2i+1}-3\right)/8\right]\left(mod\, p^{2i+1}\right)$,

\quad{}\qquad{}if $p\equiv17\left(mod\,24\right)$ $\Rightarrow$
$q\equiv9\left(mod\,24\right)$, $m_{1}\neq\left[3\left(p^{2i+1}-1\right)/8\right]\left(mod\, p^{2i+1}\right)$, 

\quad{}\qquad{}if $p\equiv19\left(mod\,24\right)$ $\Rightarrow$
$q\equiv3\left(mod\,24\right)$, $m_{1}\neq\left[\left(p^{2i+1}-3\right)/8\right]\left(mod\, p^{2i+1}\right)$.

(C3) if $p\equiv1\left(mod\,6\right)$ $\Rightarrow$ $q\equiv22\left(mod\,24\right)$,
$m_{1}\neq\left[11\left(p^{2i+1}-5\right)/12\right]\left(mod\, p^{2i+1}\right)$,

\hspace{0.75cm}if $p\equiv5\left(mod\,6\right)$ $\Rightarrow$$q\equiv14\left(mod\,24\right)$,
$m_{1}\neq\left[\left(7p^{2i+1}-5\right)/12\right]\left(mod\, p^{2i+1}\right)$

(C4.1) $q\equiv6\left(mod\,8\right)$, $m_{1}\neq2\left(mod\,3\right)$.

\medskip{}

For $\mu=11$, $A=1$, $B=1$, yielding $M\equiv11\left(mod\,12\right)$
$\Rightarrow$ $M=12m_{1}+11$:

(C1.3) if $q\equiv1\left(mod\,3\right)$ $\Rightarrow$ $q\equiv4\left(mod\,12\right)$,
$m_{1}\neq\left(3^{2i-1}-1\right)\left(mod\,3^{2i}\right)$,

\quad{}\qquad{}if $q\equiv2\left(mod\,3\right)$ $\Rightarrow$
$q\equiv8\left(mod\,12\right)$, $m_{1}\neq\left(2\times3^{2i-1}-1\right)\left(mod\,3^{2i}\right)$.

(C2) if $p\equiv\pm1\left(mod\,12\right)$ $\Rightarrow$ $q\equiv11\left(mod\,12\right)$,
$m_{1}\neq\left[11\left(p^{2i}-1\right)/12\right]\left(mod\, p^{2i}\right)$, 

\hspace{0.75cm}if $p\equiv5\left(mod\,12\right)$ $\Rightarrow$
$q\equiv7\left(mod\,12\right)$, $m_{1}\neq\left[\left(7p^{2i+1}-11\right)/12\right]\left(mod\, p^{2i+1}\right)$

\hspace{0.75cm}or $m_{1}\neq\left[\left(7\left(p-5\right)/12\right)\left(\sum_{j=0}^{2i}C_{j}^{2i+1}\left(p-5\right)^{2i-j}5^{j}\right)+\left(7\times5^{2i+1}-11\right)/12\right]$

\quad{}\qquad{}$\left(mod\, p^{2i+1}\right)$,

\hspace{0.75cm}if $p\equiv7\left(mod\,12\right)$ $\Rightarrow$
$q\equiv5\left(mod\,12\right)$, $m_{1}\neq\left[\left(5p^{2i+1}-11\right)/12\right]\left(mod\, p^{2i+1}\right)$

\hspace{0.75cm}or $m_{1}\neq\left[\left(5\left(p-7\right)/12\right)\left(\sum_{j=0}^{2i}C_{j}^{2i+1}\left(p-7\right)^{2i-j}7^{j}\right)+\left(5\times7^{2i+1}-11\right)/12\right]$

\quad{}\qquad{}$\left(mod\, p^{2i+1}\right)$.

(C3) $q\equiv0\left(mod\,12\right)$, $m_{1}\neq-1\left(mod\, p^{2i+1}\right)$.

(C4.2) $q\equiv2\left(mod\,3\right)$, $m_{1}\neq\left(3\times2^{\alpha-2}-1\right)\left(mod\,2^{\alpha}\right)$,
$\forall\alpha\geq2$. 

\medskip{}

Note that these are still only necessary conditions and that they
do not exclude all values of $M$ that do not give solutions to (\ref{eq:53-3}).
For example, as signaled by Beeckmans, the value $M=842=24\times35+2$,
although complying with the conditions (C1.3), (C2) and (C3) for $\mu=2$,
does not yield solutions to (\ref{eq:53-3}).

\section{Conclusion}

The allowed congruent values of $M$ for which sums of $M$ consecutive
integer squares equal squared integers, namely$M\equiv0,9,24$ or
$33\left(mod\,72\right)$, $M\equiv1,2$ or $16\left(mod\,24\right)$,
$M\equiv11\left(mod\,12\right)$, were further characterized with
additional congruence conditions using Beeckmans' eight necessary
conditions.

\end{document}